
\magnification1200
\input amssym.def  
\def\SetAuthorHead#1{}
\def\SetTitleHead#1{}
\def\NoindentAfter{\everypar={\setbox0=\lastbox\everypar={}}}
\def\H#1\par#2\par{{\baselineskip=15pt\parindent=0pt\parskip=0pt
 \leftskip= 0pt plus.2\hsize\rightskip=0pt plus.2\hsize
 \bf#1\unskip\break\vskip 4pt\rm#2\unskip\break\hrule
 \vskip40pt plus4pt minus4pt}\NoindentAfter}
\def\HH#1\par{{\bigbreak\noindent\bf#1\medbreak}\NoindentAfter}
\def\HHH#1\par{{\bigbreak\noindent\bf#1\unskip.\kern.4em}}
\def\th#1\par{\medbreak\noindent{\bf#1\unskip.\kern.4em}\it}
\def\endth{\medbreak\rm}
\def\pf#1\par{\medbreak\noindent{\it#1\unskip.\kern.4em}}
\def\df#1\par{\medbreak\noindent{\it#1\unskip.\kern.4em}}
\def\enddf{\medbreak}

\def\qedbox{\vrule width2mm height2mm\hglue1mm\relax}
\def\qed{\ifmmode\qedbox\else\hglue5mm\unskip\hfill\qedbox\medbreak\fi\rm}

\let\Overline\overline
\def\Smallfonts{}
\def\tenpoint{}
\def\FigureTitle#1{\medskip\centerline{\bf#1\unskip}}
\let\Item\item
\def\cite#1{{\bf[#1]}}
\def\Em#1{{\it #1\/}}
\def\Bib#1\par{\bigbreak\bgroup\centerline{#1}\medbreak\parindent30pt
 \parskip2pt\frenchspacing\par}
\def\endBib{\par\egroup}
\newdimen\Overhang
\def\rf#1{\par\noindent\hangafter1\hangindent=\parindent
     \setbox0=\hbox{[#1]}\Overhang\wd0\advance\Overhang.4em\relax
     \ifdim\Overhang>\hangindent\else\Overhang\hangindent\fi
     \hbox to \Overhang{\box0\hss}\ignorespaces}

\SetTitleHead{On irregular links at infinity of algebraic plane curves}
\SetAuthorHead{Walter D. Neumann and Le Van Thanh}

\def\C{\Bbb C}
\def\Q{\Bbb Q}
\def\P{\Bbb P}
\def\L{{\cal L}}
\def\inv{^{-1}}
\def\frame#1{\hbox{$\vcenter{\hbox{\vrule\vbox {\hrule\kern2pt%
  \hbox{\kern2pt#1\kern2pt}\kern2pt\hrule}\vrule}}$}}
\def\implies{$\Rightarrow$}
\let\O\Overline

\H On irregular links at infinity of algebraic plane curves

Walter D. Neumann and Le Van Thanh
\footnote{}{
\hfill\break
The authors thank the Max-Planck-Institut f\"ur Mathematik in Bonn,
where this work was done, for its support. The first author also
acknowledges partial support by the NSF.}

{\Smallfonts\narrower
\HHH Abstract

We give two proofs of a conjecture of Neumann \cite{N1} that a reduced
algebraic plane curve is regular at infinity if and only if its link
at infinity is a regular toral link.  This conjecture has also been
proved by Ha H.~V. \cite {H} using Lojasiewicz numbers at infinity.
Our first proof uses the polar invariant and the second proof uses
linear systems of plane curve singularities.  The second approach also
proves a stronger conjecture of \cite{N1} describing topologically the
regular link at infinity associated with an irregular link at infinity.\par}

\tenpoint
\HH 1. Introduction

Let $P:\C^2\rightarrow\C$ be a polynomial such that $V=P\inv(0)$ is a
reduced curve.  The intersection of $V$ with a sphere of sufficiently
large radius centered at the origin is transverse and gives a link
$\L(V,\infty):=(S^3,V\cap S^3)$, called the \Em{link at infinity of
$V$}.  It is independent of the radius of $S^3$ (up to isotopy) and is
a toral link.  We shall abbreviate it $\L=\L(V,\infty)=(S^3,L)$.

One says that $V$ is \Em{regular at infinity} if the defining polynomial $P$
gives a trivial fibration in a neighborhood of $V$ at infinity.  Only finitely
many fibers of $V$ are irregular at infinity.

We recall the construction in \cite{N1} of the splice diagram
$\Omega=\Omega(\L)$ for the link at infinity.  (See \cite{N1} or
\cite{EN} for basics on splice diagrams)

By compactifying we can write $V=\O V - (\O V \cap \P^1)\subset \P^2 -
\P^1 = \C^2$.  $\O V$ meets $\P^1$ in finitely many points
$Y_1,\ldots,Y_n$, say.  Let $D_0$ be a $2$-disk in $\P^1$ which
contains $\O V\cap\P^1$ and $D$ a $4$-disk neighborhood of $D_0$ in
$\P^2$ whose boundary $S=\partial D$ meets $\O V\cup \P^1$
transversally (see Fig.~1).
\midinsert
\vglue 1in
\FigureTitle{Figure 1.}
\endinsert

Then $\L_0= (S,(\P^1\cup\O V)\cap S)$ is a link which can be
represented by a splice diagram $\Gamma$ as follows:

\vglue 1in\noindent
Here $K_0$ is the component $\P^1\cap S$ and each $\longleftarrow
\kern-6pt\frame{$\Gamma_i$}$ is the diagram representing the local
link of $\P^1\cup\O V$ at the point $Y_i$.

Let $N\P^1$ be a closed tubular neighborhood of $\P^1$ in $\P^2$ whose
boundary $S^3=\partial N\P^1$ is the sphere at infinity and
$\L'=(S^3,S^3\cap V)$ is, but for orientation, the link at infinity
that interests us.  We may assume that $N\P^1$ is obtained from $D$ by
attaching a $2$-handle along $K_0$, so $\L'$ is obtained from $\L_0$
by $(+1)$-Dehn surgery on $K_0\subset S=\partial D$.

As in \cite{N1}, we call a weight in the splice diagram \Em{near} or
\Em{far} according as it is on the near or far end of its edge,
viewed from $K_0$.  As described in \cite{N1}, the splice diagram for
$\L'$ is

\vglue 1in\noindent
where $\Gamma_i'$ is obtained from $\Gamma_i$ by replacing each far
weight $\beta_v$ by $\beta_v-\lambda_v^2\alpha_v$ with
$\alpha_v$ equal to the product of the near weights at vertex $v$ and
$\lambda_v$ the product of the weights adjacent to, but not on, the
simple path from $v$ to the vertex corresponding to $K_0$.

Finally, we must reverse orientation to consider $S^3$ as a large
sphere in $\C^2$ rather than as $\partial N\P^1$. The effect is to
reverse the signs of all near weights.  We can then forget the
leftmost vertex, which is redundant, to get a diagram $\Omega$ for
$\L=\L(V,\infty)$ as follows.
\vglue 1in
\noindent The leftmost vertex of $\Omega$ is called the \Em{root
vertex}.  $\Omega$ is an \Em{RPI splice diagram} in the sense of \cite{N1}.
``RPI'' stands for ``reverse Puiseux inequalities'' and refers to certain
inequalities that the weights of $\Omega$ satisfy.

For each vertex $v$ of $\Omega$ we denote by $\delta_v$ its valency
(number of incident edges).  We have three types of vertices:
\Em{arrowheads} (corresponding to components of $\L$), \Em{leaves}
(non-arrowheads of valency $1$) and \Em{nodes} (non-arrowheads of
valency $>2$).  For each non-arrowhead, denote by $S_v$ a
corresponding \Em{virtual component} for $\L$ (see \cite{N1}).  In
particular, $S_o=K_0$, where $o$ denotes the root vertex.  The linking
number $b_v := l(S_v,K_0)$ is called the \Em{braid index} of $v$ and
$l_v := l(S_v,L)$ (sum of linking numbers of $S_v$ with all components
of $\L=(S^3,L)$) is the \Em{(total) linking coefficient} at $v$
(called ``multiplicity at $v$ in \cite{EN}). The linking coefficient
$l_o$ at the root vertex is the degree $d$ of the defining polynomial.

\df Definition 1.1

$\Omega$ is a \Em{regular} RPI splice diagram if $l_v\ge 0$ for all
non-arrowhead vertices.

As described in \cite{N1}, the diagram $\Omega$ depends of the
particular compactification $\C^2\subset\P^2$ of $\C^2$ that we
choose, but if one RPI splice diagram for $\L$ is regular then every
RPI splice diagram for $\L$ is.  Thus the regularity or irregularity
of $\Omega$ is a topological property of $\L$, and we say the toral
link $\L$ is \Em{regular} or \Em{irregular} accordingly.
\enddf

In \cite{N1} it was shown that $\L(V,\infty)$ is a regular link if $V$
is regular at infinity, and the converse was conjectured.  This conjecture
has been proved by Ha~H.~V. in \cite{H}.
One purpose of this note is to give a new short proof of this using
completely different methods, namely the polar invariant of $V$.  We shall
prove

\th Theorem 1.2

The following conditions are equivalent:
\Item{(1)} $V=P\inv (0)$ is regular at infinity.
\Item{(2)} $\L(V,\infty)$ is a regular toral link.
\Item{(3)} For every component $\gamma$ of the polar curve for $V$,
the natural valuation at infinity on holomorphic functions on $\gamma$
satisfies
$v_\gamma^\infty(P|\gamma)\ge 0$.
\Item{(4)} There exist $R>0$ and $\delta>0$ and a line
field $\alpha{\partial\over\partial x} + \beta{\partial\over\partial
y}$ which is transverse to $P\inv(t)\cap(\C^2-B_R)$ for all
$t<\delta$, where $B_R=\{(x,y)\in\C^2:|(x,y)|<R\}$.
\endth

Using different methods we shall also describe the proof of a stronger
conjecture from \cite{N1}, which computes the link at infinity of any
regular fiber of the defining polynomial $P$ in terms of any RPI
splice diagram $\Omega$ for an irregular link at infinity of $P$.
Namely, let $\Omega^-$ be the full subgraph of $\Omega$ on all
vertices $v$ with $l_v<0$ and arrowheads adjacent to them, and let
$\Omega^-_1,\ldots,\Omega^-_k$ be the connected components of
$\Omega^-$.  For each $j=1,\ldots,k$, $\Omega^-_j$ is connected to the
rest of $\Omega$ by a single edge.  It is easy to see that there is a
unique way of cutting at this edge and replacing $\Omega^-_j$ by a
graph of the form
\vglue 1in
\noindent in such a way that in the new $\Omega$:
\Item{(i)} $l_w=0$, where $w$ is the new node;
\Item{(ii)} $l_v$ is unchanged for every non-arrowhead vertex of
$\Omega-\Omega^-_j$.

\noindent Let $\Omega_0$ be the result of doing this for each $j=1,\ldots,k$.

\th Theorem 1.3

$\Omega_0$ is an RPI splice diagram for the link at infinity of any
regular fiber of $P$.\endth

\HH 2.~~Proof of Theorem 1.2

That (1)\implies(2) was proved in \cite{N1}.  That (4)\implies(1) is
trivial: if (4) holds then the vector field
$(\alpha{\partial\over\partial x}+\beta{\partial\over\partial y})/
(\alpha{\partial P\over\partial x}+\beta{\partial P\over\partial y})$
trivializes the neighborhood at infinity $P\inv
(\Delta_{\delta})\cap(\C^2-B_R)$ of $V=P\inv (0)$, where
$\Delta_{\delta}=\{t\in\C:|t|<\delta\}$ and
$B_R=\{(x,y)\in\C^2:|(x,y)|<R\}$.

Suppose (3) holds, that is, for all components $\gamma$ of the polar
curve $$\Pi=\{(x,y):\alpha{\partial P\over\partial x}+\beta{\partial
P\over\partial y}=0\}\quad\hbox{$\alpha$ and $\beta$ sufficiently
general},$$ one has $v_\gamma^\infty(P|\gamma)\ge 0$.  This means that
for all $\gamma\subset\Pi$, $P(x,y)|\gamma\not\rightarrow0$ as
$|(x,y)|\rightarrow\infty$, so there exists $\delta>0$ such that for
sufficiently large $R$ one has:
$$P\inv(\Delta_{\delta})\cap(\C^2-B_R)\cap\Pi=\emptyset,$$
This is statement (4) of Theorem 1.2.

It remains to prove that (2)\implies(3).  Recall that $\Gamma_i$ is
the splice diagram for the local link $\L_i=(S^3_\epsilon(Y_i),
S^3_\epsilon(Y_i)\cap\O V)$ of $\O V\subset\P^2$ at the point $Y_i$.
If $v$ is a non-arrowhead vertex of $\Gamma_i$, let $S_v$ be the
corresponding virtual component of $\L_i$. Write $\L_i=(S^3,L_i)$, and
denote the ``local braid index'' and ``local linking coefficient'' at
$v$ by
$$\eqalign{ b_v^{loc}&:=l_{\L_i}(S_v,K_0),\cr
l_v^{loc}&:=l_{\L_i}(S_v,L_i).\cr} $$
(Note that this local braid index is the braid index with respect to
the special line $z=0$ rather than with respect to a general line,
which is how the term is usually used.)  Since we can consider $v$ as
a vertex of $\Omega$, the (global) braid index $b_v$ and linking
coefficient $l_v$ are also defined.

\th Lemma 2.1

\hfill\llap{\hbox to \hsize{\hskip 1 in{\rm (i)}~~$b_v=b_v^{loc}$\hfil}}

\noindent\hskip 1in{\rm(ii)}~~$l_v=b_vd-l_v^{loc}$.
\endth

\pf Proof

Recall (e.g., Lemma 3.2 of \cite{N1}) that the linking number of any
two components (virtual or genuine) of a toral link is the product of
all weights adjacent to but not on the simple path joining the
corresponding vertices of a splice diagram for the link.  Given the
relationship between the weights of $\Gamma$ and $\Omega$ described in
Section 1, (i) is  then immediate and (ii) is a simple calculation.

Both (i) and (ii) are also easy to see topologically.  For (ii) one
uses the fact that $(1/d)$-Dehn surgery on a knot $K$ in $S^3$
replaces the linking number $l(C_1,C_2)$ of any two disjoint
$1$-cycles which are disjoint from $K$ by
$l(C_1,C_2)-dl(C_1,K)l(C_2,K)$.  Applying this with $C_1=S_v$,
$C_2=L$, we see, in the notation of Section 1, $l_v=l_\L(S_v,L) =
-l_{\L'}(S_v,L)=-(l_{\L_0}(S_v,L)-l(S_v,K_0)l(L,K_0)) =
-(l_v^{loc}-b_vd)$, since
$l_{\L_0}(S_v,L)=l_{\L_0}(S_v,L_i)=l_{\L_i}(S_v,L_i)=l_v^{loc}$.\qed

Let us assume part (3) of Theorem 1.2 is not true.  We may assume
coordinates are chosen such that the polar curve is given as
$$\Pi=\{(x,y)\in\C^2:{\partial P(x,y)\over\partial y}=0\}$$
and, in addition, the component $\gamma\subset\Pi$ contravening (3)
has Puiseux expansion at infinity
$$\gamma=\{(x,y):y=y_\gamma(x)\}\quad\hbox{with}\quad
y_\gamma(x)=\sum_{\alpha\in1+\Q^-}a_\alpha x^\alpha,$$
where $\Q^-=\{\alpha\in\Q:\alpha<0\}$ and $P_y(x,y_\gamma(x))\equiv
0$.  In particular, the point at infinity in question is $y=0$.

One has
$$v_\gamma^\infty(P|\gamma)=v_x^\infty(P(x,y_\gamma(x))).$$
On the other hand, if $\tilde P(z,y)=z^dP({1\over z},{y\over z})$
is the equation for $\Overline{P^{-1}(0)}$ in local coordinates at infinity,
then
$\tilde P_y(z,zy_\gamma({1\over  z}))\equiv 0$,
which means that
$$\tilde\gamma := \{(z,y):y=zy_\gamma({1\over z})=:
y_{\tilde\gamma}(z)\}$$
is the compactification of $\gamma$.  One has
$v_z(\tilde P(z,y_{\tilde\gamma}(z))=d-v_x^\infty(P(x,y_\gamma(x))$,
so, by assumption,
$$v_z(\tilde P(z,y_{\tilde\gamma}(z))>d.$$
But, in the local situation, one knows that
$$v_z(\tilde P(z,y_{\tilde\gamma}(z)))={\tilde\gamma\cdot\tilde P\inv
(0)\over \tilde\gamma\cdot\{z=0\}}.$$

If $\{z=0\}$ is a generic line, then the main result of \cite{LMW} can
be expressed that for some vertex $v$ of the splice diagram $\Gamma$
(which is the diagram for the link at $0$ of $z\tilde P(z,y)=0$, rather than
$\tilde P(z,y)=0$ as in \cite{LMW}),
$${\tilde\gamma\cdot\tilde P\inv(0)\over
\tilde\gamma\cdot\{z=0\}}= {l_v^{loc}\over b_v^{loc}}.$$
Expressed in this form, the proof in \cite{LMW} applies also if
$\{z=0\}$ is not a generic line
(cf.~\cite{L}).
Thus $l_v^{loc}/b_v^{loc}>d$, so $\L$ is not
regular by Lemma 2.1 (ii).\qed

\HH 3.~~Proof of Theorem 1.3

Let $f(y,z)$ be a polynomial with $f(0,0)=0$ and $a$ and $b$
non-negative integers and consider the linear family of polynomials
$f_t(y,z)=f(y,z)-ty^az^b$ with $t\in\C$.  Let $V_t= f_t\inv(0)$. A
special case of the main theorem of \cite{N2} describes the local link
at $0$ of a generic member $V_t$ of this family in terms of the link
at $0$ of $V_0$.  Namely, we can assume without loss of generality
that $f(y,z)$ is not divisible by $y$ or $z$.  Let $\Gamma$ be the
splice diagram for the link of $yzf(y,z)=0$ at $0$ and let $K_y$ and
$K_z$ be the components of this link given by $y=0$ and $z=0$ and $K$
the union of the remaining components (that is, the link of $f=0$).
For each vertex $v$ of $\Gamma$ consider the number $s_v :=
al(S_v,K_y)+bl(S_v,K_z)-l_v$.

Let $\Gamma^-$ be the full subgraph of $\Gamma$ on all vertices $v$
with $s_v<0$ and arrowheads adjacent to them, and let
$\Gamma^-_1,\ldots,\Gamma^-_k$ be the connected components of
$\Gamma^-$.  Suppose:
\Item{(*)} For each $j=1,\ldots,k$, $\Gamma^-_j$ is
connected to the rest of $\Gamma$ by a single edge.

\noindent
In \cite{N2} it is shown that here is a unique way of cutting at this
edge and replacing $\Gamma^-_j$ by a graph of the form
\vglue 1in
\noindent in such a way that in the new $\Gamma$:
\Item{(i)} $s_w=0$, where $w$ is the new node;
\Item{(ii)} $s_v$ is unchanged for every non-arrowhead vertex of
$\Gamma-\Gamma^-_j$.

\noindent
Moreover, the result $\Gamma_0$ of doing this for each $j=1,\ldots,k$
is the splice diagram for the link at 0 of $yzf_t(y,z)=0$ for generic
$t$.
(If (*) fails then $\Gamma^-$ is connected and meets the rest of $\Gamma$
in two edges and $\Gamma_0$ is computed similarly in \cite{N2}.)

Now let $\tilde P(z,y)=0$ be the polynomial at a point at infinity for
a curve $P(x,y)=0$, as described in the previous section. Put
$f(y,z)=\tilde P(z,y)$, $a=0$, $b=\deg(P)=d$.  Then $f_t(y,z)=0$ is
the equation of the point at infinity for $P(x,y)=t$.  By Lemma
2.1(ii), and the relationship between the splice diagrams $\Gamma$ and
$\Omega$ described in section 2, the above result of \cite{N2}
translates directly to give Theorem 1.3.

\Bib Bibliography

\rf{EN}  D. Eisenbud and W. D. Neumann,  {\it Three-Dimensional Link Theory and
Invariants of Plane Curve Singularities}, Annals of Math. Studies {\bf
110} (Princeton Univ. Press 1985).

\rf{H} Ha H.V., On the irregular at infinity algebraic plane curve, (Preprint
91/4, Institute of Math., National Center for Scientific Research of Vietnam).

\rf{L} L. V. Than, Affine polar quotients of algebraic plane curves, (in
preparation).

\rf{LMW} Le D.T., F. Michel, C. Weber, Courbes polaires et topologie des
courbes planes, Ann. Sc. Ec. Norm. Sup. $4^e$ Series, {\bf 24} (1991),
141--169.

\rf{N1} W. D. Neumann, Complex plane curves via their links at infinity,
Invent.
Math. {\bf 98} (1989), 445--489.

\rf{N2} W. D. Neumann, Linear systems of plane curve singularitites and
irregular
links at infinity, (in preparation).

\endBib

\bigskip
\indent\vbox{\hsize .4\hsize\parindent=0pt\Smallfonts \obeylines
Walter D. Neumann
Department of Mathematics
Ohio State University
Columbus, OH 43210, USA
neumann@mps.ohio-state.edu
}
\vbox{\hsize .4\hsize\parindent=0pt\Smallfonts \obeylines
Le Van Than
Institute of Mathematics
P.O. Box 631-10000
Hanoi, Vietnam
{}~~
}
\bye